\documentclass[10pt]{article}

\usepackage{amsmath}
\usepackage{amssymb}
\usepackage{indentfirst}
\usepackage{graphics} 
\usepackage{color}

\makeatletter

\@addtoreset{equation}{section}
\makeatother
\def\<{\langle}
\def\>{\rangle}

\newtheorem{lem}{Lemma}[section]
\newtheorem{theo}{Theorem}[section]
\newtheorem{rem}{Remark}[section]

\makeatletter
   
   \@addtoreset{equation}{section}
\makeatother

\begin{document}
\title{\bf A role of potential on $L^{2}$-estimates\\ for some evolution equations}
\author{Ryo Ikehata\thanks{Corresponding author: ikehatar@hiroshima-u.ac.jp} \\ {\small Department of Mathematics, Division of Educational Sciences}\\ {\small Graduate School of Humanities and Social Sciences} \\ {\small Hiroshima University} \\ {\small Higashi-Hiroshima 739-8524, Japan} 
}
\date{}
\maketitle
\begin{abstract}
In this papwe we consider an effective role of the potential of the wave equations with/without damping on the $L^{2}$-estimate of the solution itself. In the free wave equation case it is known (\cite{ike-2023}) that the $L^{2}$-norm of the solution itself generally grows to infinity (as $t \to \infty$) in the one and two dimensional cases, however, by adding the potential with quite generous condition one can controle the growth property to get the $L^{2}$-bounds. This idea can be also applied to the damped wave equations with potential in order to get fast energy and $L^{2}$ decay results in the low dimensional case, which are open for a long period. Applications to heat and plate equations with a potential can be also studied. In this paper the low dimensional case is a main target.
\end{abstract}
\section{Introduction}
\footnote[0]{Keywords and Phrases: Wave, heat ad plate equations; potential; Cauchy problem, $L^{2}$-bound, damping, fast decay.}
\footnote[0]{2010 Mathematics Subject Classification. Primary 35L05; Secondary 35B40, 35L30.}

In this paper, first we are concerned with the following Cauchy problem with a potential:

\begin{equation}\label{01}
u_{tt}- \Delta u + V(x)u = 0, \quad t > 0,\quad x \in {\bf R}^{n},
\end{equation}
\begin{equation}\label{02}
u(0,x)=u_{0}(x),\hspace{0,50cm}u_{t}(0,x) = u_{1}(x)\quad x \in {\bf R}^{n}.
\end{equation} 
Throughout the paper one assumes that (for simplicity) the potential function $V \in BC({\bf R}^{n})$ satisfies $V (x) \geq 0$ ($x \in {\bf R}^{n}$).\\ 
Note that functions and solutions treated in this paper are all real-valued.\\
\noindent
One defines the total energy for the equation \eqref{01} by
\begin{equation}\label{E1}
E(t) := \frac{1}{2}\int_{{\bf R}^{n}}\left(\vert u_{t}(t,x)\vert^{2} + \vert\nabla u(t,x)\vert^{2} + V(x)\vert u(t,x)\vert^{2} \right)dx.
\end{equation}
Then, it is known that for each initial data $[u_{0},u_{1}] \in H^{1}({\bf R}^{n}) \times L^{2}({\bf R}^{n})$ the problem \eqref{01}-\eqref{02} has a unique weak solution $u \in C([0,\infty); H^{1}({\bf R}^{n})) \cap C^{1}([0,\infty); L^{2}({\bf R}^{n}))$ satisfying the energy identity 
$$E(t) = E(0).$$

In the first part of this paper, we consider the $L^{2}$-boundedness in time of the solution itself to problem \eqref{01}-\eqref{02}. Why $L^{2}$-boundedness are concerned is that as is pointed out recently in \cite{ike-2023}, in the case of $\int_{{\bf R}^{n}}u_{1}(x)dx \ne 0$, the solution itself to the free wave equation grows up at a suitable rate in the case of $n = 1,2$:
\begin{equation}\label{04-1}
\Vert w(t,\cdot)\Vert \approx \sqrt{t}\, (n = 1),\quad \Vert w(t,\cdot)\Vert \approx \sqrt{\log t}\, (n = 2),\,\,\,(t \gg 1),
\end{equation}
where $w(t,x)$ is the solution to the free waves with initial data $[u_{0},u_{1}] \in C_{0}^{\infty}({\bf R}^{n}) \times C_{0}^{\infty}({\bf R}^{n})$ (for simplicity):
\begin{equation}\label{03}
w_{tt}- \Delta w = 0, \quad t > 0,\quad x \in {\bf R}^{n},
\end{equation}
\begin{equation}\label{04}
w(0,x)=u_{0}(x),\hspace{0,50cm}w_{t}(0,x) = u_{1}(x)\quad x \in {\bf R}^{n}.
\end{equation} 
So, it seems not trivial to get even $L^{2}$-bounds of the solution itself to problem \eqref{01}-\eqref{02} if we are treating a wealky effective potential case, for example, $V(x) = e^{-\vert x\vert^{2}}$, and so on. Also, as is studied in previous researchs (\cite{B}, \cite {D}, \cite{GV, GV-2}, \cite{pe}, \cite{P} and the references therein) it seems important and difficult to get $L^{p}$-bounds and/or $L^{p}$-decay estimates of the solutions. In fact, it is known in \cite{GV} that for $n = 3$, to get $L^{4}$-estimate of the solution one must restrict the decay condition (as $\vert x\vert \to \infty$) assumed on the potential in a context. In such cases, one can not treat low dimensions and smaller effective potential even for $L^{2}$-bound. In general, as compared to the higher dimensional case ($n \geq 3$), it seems that the low dimensional case $n = 1,2$) is less studied. In the former part, by restricting the size of initial data one can controle such growth property by adding the potential. The main result of this paper reads as follows.
\begin{theo}\label{t01} 
Let $n\geq 1$, and $V \in BC({\bf R}^{n})$ satifies $V(x) > 0$ {\rm (}$\forall x \in {\bf R}^{n}${\rm )}. Suppose that the initial data $(u_{0}, u_{1}) \in H^{1}({\bf R}^{n})\times L^{2}({\bf R}^{n})$ further satisfies
\[I_{0}^{2} := \int_{{\bf R}^{n}}\frac{\vert u_1(x)\vert^{2}}{V(x)}dx < +\infty.\]
Then the unique solution $u \in C([0,\infty); H^{1}({\bf R}^{n})) \cap C^{1}([0,\infty); L^{2}({\bf R}^{n}))$ to problem \eqref{01}-\eqref{02} satisfies
\[\Vert u(t,\cdot)\Vert  \leq C(\Vert u_{0}\Vert + I_{0}), \quad t \geq 0,\]
with some generous constant $C > 0$.
\end{theo}
{\bf Example 1.}{\rm The constant coefficient case $V(x) = m^{2}$ ($m > 0$) can be included as a trivial example since $I_{0} < +\infty$ is automatically satisfied. This is the so-called Klein-Gordon equation case.}\\
\noindent
{\bf Example 2.} {\rm One can also choose $V(x) := e^{-\vert x\vert^{2}}$ as the 2nd example. In this case one has to assume $I_{0}^{2} = \int_{{\bf R}^{n}}e^{\vert x\vert^{2}}\vert u_{1}(x)\vert^{2}dx < +\infty$ additionally to get the $L^{2}$-boundedness. The potential $V(x) := e^{-\vert x\vert^{2}}$ can be viewed as $V(x) \approx 0$ for $\vert x\vert \geq L$ for large $L > 0$.}\\
{\bf Example 3.} {\rm As the 3rd example, one can choose $V(x) := (1+\vert x\vert^{2})^{-\alpha/2}$ ($\alpha > 0$). In this case one has to assume $I_{0}^{2} = \int_{{\bf R}^{n}}(1+\vert x\vert^{2})^{\alpha/2}\vert u_{1}(x)\vert^{2}dx < +\infty$ additionally to get the $L^{2}$-boundedness.}\\

So, one can ask a question: can one localize the effective region of $V(x)$ to the smaller region to still have $L^{2}$-boundedness. A next result gives an answer for this question. 

\begin{theo}\label{t02} 
Let $n\geq 1$. Let $\Omega \subset {\bf R}^{n}$ be an open set, and we assume that $V \in BC({\bf R}^{n})$ satifies $V(x) > 0$ for $\forall x \in \Omega$, and $V(x) \geq 0$ for $\forall x \in {\bf R}^{n}$. If $(u_{0}, u_{1})\in H^{1}({\bf R}^{n})\times L^{2}({\bf R}^{n})$ additionally satisfies ${\rm supp}u_{1} \subset \Omega$ and
\[J_{0}^{2} := \int_{\Omega}\frac{\vert u_1(x)\vert^{2}}{V(x)}dx < +\infty,\]
then the unique weak solution $u \in C([0,\infty); H^{1}({\bf R}^{n})) \cap C^{1}([0,\infty); L^{2}({\bf R}^{n}))$ to problem \eqref{01}-\eqref{02} satisfies
\[\Vert u(t,\cdot)\Vert \leq C(\Vert u_{0}\Vert + J_{0}),\]
with some generous constant $C > 0$.
\end{theo}
{\bf Example 4.} {\rm If $\Omega = B_{\varepsilon}$ with so small small $\varepsilon > 0$, under the assumption that ${\rm supp}u_{1} \subset B_{\varepsilon}$, and $\int_{B_{\varepsilon}}\frac{\vert u_1(x)\vert^{2}}{V(x)}dx < +\infty$, then one can still get the $L^{2}$-boundedness.}\\

As pointed out in Example 4, even if the effective region of the potential $V(x)$ is small, one can get the $L^{2}$-bounds of the solution $u(t,x)$ to \eqref{01}-\eqref{02} and growing up never occurs, that is, there exists a quite big difference between the free waves and the wave equation with a localized potential from the $L^{2}$-behavior point of view.

In the second part of this paper, one can also apply the method developed in the former part to the damped wave equation with a smaller potential:
\begin{equation}\label{01-w}
u_{tt}- \Delta u + V(x)u + u_{t} = 0, \quad t > 0,\quad x \in {\bf R}^{n},
\end{equation}
\begin{equation}\label{02-w}
u(0,x)=u_{0}(x),\hspace{0,50cm}u_{t}(0,x) = u_{1}(x)\quad x \in {\bf R}^{n}.
\end{equation} 
\noindent
By applying the same concept as Theorems \ref{t01} and \ref{t02} one can get the following results.
\begin{theo}\label{t11} 
Let $n\geq 1$, and $V \in BC({\bf R}^{n})$ satisfies $V(x) > 0$ {\rm (}$\forall x \in {\bf R}^{n}${\rm )}. Suppose that the initial data $(u_{0}, u_{1}) \in H^{1}({\bf R}^{n})\times L^{2}({\bf R}^{n})$ further satisfies
\[K_{0}^{2} := \int_{{\bf R}^{n}}\frac{\vert u_1(x)+u_{0}(x)\vert^{2}}{V(x)}dx < +\infty.\]
Then the unique solution $u \in C([0,\infty); H^{1}({\bf R}^{n})) \cap C^{1}([0,\infty); L^{2}({\bf R}^{n}))$ to problem \eqref{01-w}-\eqref{02-w} satisfies
\[(1+t)^{2}E(t) \leq C(\Vert u_{0}\Vert^{2} + \Vert u_{1}\Vert^{2} + \Vert \nabla u_{0}\Vert^{2} + K_{0}^{2}), \quad (t \geq 0),\]
\[(1+t)\Vert u(t,\cdot)\Vert^{2}  \leq C(\Vert u_{0}\Vert^{2} + \Vert u_{1}\Vert^{2} + \Vert \nabla u_{0}\Vert^{2} + K_{0}^{2}), \quad (t \geq 0)\]
with some generous constant $C > 0$.
\end{theo}
\begin{rem}\, {\rm It should be mentioned that our result can be soon applied to the more general damping coefficient case such that
\[u_{tt}- \Delta u + V(x)u + a(x)u_{t} = 0,\] 
where $a \in C({\bf R}^{n})$ satisfies
\[0 < a_{0} \leq a(x) \leq b_{0} < +\infty\quad (x \in {\bf R}^{n})\]
and $a_{0}$ and $b_{0}$ are constants. This is one of our merits. In this case, the obtained results are slightly modified with several constants depending on such $a_{0}$ and $b_{0}$. In connection with this, an application of our method to the localized damping case studied in \cite{ike-2003} will be a forthcoming project.}
\end{rem}
\begin{rem}
{\rm Note that
\[\int_{{\bf R}^{n}}\vert u_{1}(x)+u_{0}(x)\vert dx \leq K_{0}(\int_{{\bf R}^{n}}V(x)dx)^{1/2}.\]
Thus, if $V \in L^{1}({\bf R}^{n})$, then one has $u_{1}+u_{0} \in L^{1}({\bf R}^{n})$. In this case, one can compare the Matsumura result \cite{m} in the $L^{1}$-framework. In the $n = 2$ case, the obtained result in Theorem \ref{t11} seems optimal, however, $n = 1$ case, surprisingly one can get faster decay estimate rather than the corresponding Matsumura estimate with $V(x) \equiv 0$. This may show a stronger effect of the $L^{1}$-potential function $V(x)$ can be captured. It should be emphasized that our hypothesis assumed on the initial data $u_{1}+u_{0}$ does not necessarily imply $u_{1} + u_{0} \in L^{1}({\bf R}^{n})$ even if one gives the additional condition $K_{0} < +\infty$ in response to the potential $V(x)$.}
\end{rem}
\begin{rem}
{\rm One can compare our results with the ones in \cite{N}. Nakao \cite{N} treats the potential $V(x) := k_{3}(1+\vert x\vert)^{-\theta}$ ($k_{3} > 0$), and the author in \cite[(iii) of Theorem 2.1]{N} obtains the decay estimate
\[0 \leq \theta \leq 1 \Rightarrow E(t) \leq C(E(0))e^{-k't^{1-\theta}}\]
for some $k' > 0$. This shows the exponetial decay, however, if one choose $\theta = 1$ one just gets the bounded estimate on $t$, while, if one applies Theorem \ref{t11}, one can get the algebraic decay for all $\theta \geq 0$ and all $n \geq 1$. Additionally, our augument is free from the finite speed of propagation property as was used in \cite{N}. One can also apply the result obtained in \cite[Theorem 1]{N-2} to get exponential total energy decay, however, at that case only $\theta \in [0,1/2)$ can be treated.}
\end{rem}

One can get the following result similarly to Theorem 1.2. The effect of potential is essential in the result.
\begin{theo}\label{t12} 
Let $n\geq 1$. Let $\Omega \subset {\bf R}^{n}$ be an open set, and we assume that $V \in BC({\bf R}^{n})$ satisfies $V(x) > 0$ for $\forall x \in \Omega$, and $V(x) \geq 0$ for $\forall x \in {\bf R}^{n}$. If $(u_{0}, u_{1})\in H^{1}({\bf R}^{n})\times L^{2}({\bf R}^{n})$ additionally satisfies ${\rm supp}\left(u_{1} + u_{0}\right) \subset \bar{\Omega}$ and
\[L_{0}^{2} := \int_{\Omega}\frac{\vert u_1(x)+u_{0}(x)\vert^{2}}{V(x)}dx < +\infty,\]
then the unique weak solution $u \in C([0,\infty); H^{1}({\bf R}^{n})) \cap C^{1}([0,\infty); L^{2}({\bf R}^{n}))$ to problem \eqref{01-w}-\eqref{02-w} satisfies
\[(1+t)^{2}E(t) \leq C(\Vert u_{0}\Vert^{2} + \Vert u_{1}\Vert^{2} + \Vert \nabla u_{0}\Vert^{2} + L_{0}^{2}), \quad (t \geq 0),\]
\[(1+t)\Vert u(t,\cdot)\Vert^{2}  \leq C(\Vert u_{0}\Vert^{2} + \Vert u_{1}\Vert^{2} + \Vert \nabla u_{0}\Vert^{2} + L_{0}^{2}), \quad (t \geq 0)\]
with some generous constant $C > 0$.
\end{theo}

\begin{rem}{\rm In \cite{i-s} and \cite{ike-2003} so-called fast decay as in Theorem \ref{t11} has been already derived for damped wave equation in an exterior domain included in ${\bf R}^{n}$ with $n \geq 2$. The restriction $n \geq 2$ comes from the use of the Hardy inequality. So the results in \cite{i-s} and \cite{ike-2003} can be also correct in the whole space case on ${\bf R}^{n}$, however, in this case the dimension $n$ must be larger than $3$ because of the use of Hardy inequality. Note that in the two dimensional case, such Hardy type inequality never holds as was pointed out by \cite{MOW}. So, it is highly desirable to investigate whether the fast decay holds or not even in the low dimensional case of the Cauchy problem to the equation with damping $a(x) \geq 0$:  
\[u_{tt}-\Delta u + a(x)u_{t} = 0, \quad t > 0, \quad x \in {\bf R}^{n}.\]
One can obtain such fast decay results by borrowing the help of potential $V(x)$ as stated in Theorem \ref{t11} and \ref{t12}. Additionally, as in Theorem \ref{t12} even if the effect of the potential is localized, one can get the fast decay. This observation seems new. If the potential $V(x)$ disappears completely, the effect of growth as in \eqref{04-1} may come into play to disturb the decay property.}
\end{rem}

At the third part of this paper, as a byproduct we can apply the method used in proving Theorem \ref{t11} to the following parabolic equations with a potential (see \cite{B} for its related problems). 
\begin{equation}\label{01-h}
u_{t}- \Delta u + V(x)u = 0, \quad t > 0,\quad x \in {\bf R}^{n},
\end{equation}
\begin{equation}\label{02-h}
u(0,x)=u_{0}(x),\quad x \in {\bf R}^{n}.
\end{equation} 
Note that in our case the operator $H := -\Delta + V(\cdot)$ is non-negative and self-adjoint in $L^{2}({\bf R}^{n})$ with its domain $D(H) = H^{2}({\bf R}^{n})$ because of the Kato-Rellich Theorem.

\begin{theo}\label{h11} 
Let $n\geq 1$, and $V \in BC({\bf R}^{n})$ satisfies $V(x) > 0$ {\rm (}$\forall x \in {\bf R}^{n}${\rm )}. Suppose that the initial data $u_{0} \in L^{2}({\bf R}^{n})$ further satisfies
\[K_{0,h}^{2} := \int_{{\bf R}^{n}}\frac{\vert u_0(x)\vert^{2}}{V(x)}dx < +\infty.\]
Then the unique solution $u \in C([0,\infty); L^{2}({\bf R}^{n})) \cap C^{1}((0,\infty); L^{2}({\bf R}^{n}))\cap C((0,\infty); H^{2}({\bf R}^{n}))$ to problem \eqref{01-h}-\eqref{02-h} satisfies
\[(1+t)\Vert u(t,\cdot)\Vert^{2}  \leq C(\Vert u_{0}\Vert^{2} + K_{0,h}^{2}), \quad (t \geq 0)\]
with some generous constant $C > 0$.
\end{theo}
\begin{rem}{\rm In \cite[Theorem 1.6]{D} the authors derive the parabolic dispersive estimate
\begin{equation}\label{h-1}
\Vert u(t,\cdot)\Vert_{q} \leq C t^{\frac{n}{2}(\frac{1}{q}-\frac{1}{p})}\Vert u_{0}\Vert_{p}, \quad \frac{1}{p} + \frac{1}{q} = 1, \quad q \in [2,\infty].
\end{equation}
In this case, the potential $V(x)$ must be the so-called of Kato class, and have a finite Kato norm (for its definition, see \cite{D}). If we choose $q = 2$ in \eqref{h-1}, then it only shows the $L^{2}$-boundedness in time. In this sense, $q = 2$ is critical in \eqref{h-1}, and seems rather delicate to be treated. So, it seems unknown whether $L^{2}$-decay property holds or not in the case of $q = 2$ for a quite general class of potential. The result of Theorem \ref{h11} says that if we choose a restricted initial data such that $K_{0,h} < +\infty$ then one can get the $L^{2}$-decay under the positivity assumption on the potential. This observation seems new. It should be emphasized that Theorem \ref{h11} holds for all $n \geq 1$.}
\end{rem}
One can also get the following result as in Theorem \ref{t12}.

\begin{theo}\label{h12} 
Let $n\geq 1$. Let $\Omega \subset {\bf R}^{n}$ be an open set, and we assume that $V \in BC({\bf R}^{n})$ satisfies $V(x) > 0$ for $\forall x \in \Omega$, and $V(x) \geq 0$ for $\forall x \in {\bf R}^{n}$. If $u_{0} \in L^{2}({\bf R}^{n})$ additionally satisfies ${\rm supp}u_{0} \subset \Omega$ and
\[J_{0,h}^{2} := \int_{\Omega}\frac{\vert u_{0}(x)\vert^{2}}{V(x)}dx < +\infty,\]
then the unique solution $u \in C([0,\infty); L^{2}({\bf R}^{n})) \cap C^{1}((0,\infty); L^{2}({\bf R}^{n}))\cap C((0,\infty); H^{2}({\bf R}^{n}))$ to problem \eqref{01-h}-\eqref{02-h} satisfies
\[(1+t)\Vert u(t,\cdot)\Vert^{2}  \leq C(\Vert u_{0}\Vert^{2} + J_{0,h}^{2}), \quad (t \geq 0)\]
with some generous constant $C > 0$.
\end{theo}

In the final part of this paper, as one more byproduct one can apply the method to the following plate equations with a potential. 
\begin{equation}\label{01-p}
u_{tt} + \Delta^{2} u + V(x)u = 0, \quad t > 0,\quad x \in {\bf R}^{n},
\end{equation}
\begin{equation}\label{02-p}
u(0,x)=u_{0}(x),\quad u_{t}(0,x) =  u_{1}(x), \quad x \in {\bf R}^{n}.
\end{equation} 
Let us first consider the known results for \eqref{01-p} in the case of $V(x) \equiv 0$. Then, in \cite{ike-2023p} the following $L^{2}$-growth estimates are known:
\begin{equation}\label{04-p}
\Vert w(t,\cdot)\Vert \approx t^{1-\frac{n}{4}}\, (n = 1,2,3),\quad \Vert w(t,\cdot)\Vert \approx \sqrt{\log t}\, (n = 4),\,\,\,(t \gg 1),
\end{equation}
where one must assume $\int_{{\bf R}^{n}}u_{1}(x)dx \ne 0$ to get \eqref{04-p}, and $w(t,x)$ is the solution to the free plate equation with initial data $[u_{0},u_{1}] \in C_{0}^{\infty}({\bf R}^{n}) \times C_{0}^{\infty}({\bf R}^{n})$ (for simplicity):
\begin{equation}\label{p03}
w_{tt} + \Delta^{2} w = 0, \quad t > 0,\quad x \in {\bf R}^{n},
\end{equation}
\begin{equation}\label{p04}
w(0,x) = u_{0}(x),\hspace{0,50cm}w_{t}(0,x) = u_{1}(x)\quad x \in {\bf R}^{n}.
\end{equation} 
So, it is natural to ask for what kind of potential comes into play for bounded estimates in time of the quantity $\Vert u(t,\cdot)\Vert$. Furthermore, in \cite[(2.2) of Theorem 2.1]{L} the author derives the estimate:
\begin{equation}\label{p-1}
\Vert u(t,\cdot)\Vert_{q} \leq C t^{\frac{n}{2q}-\frac{n}{4}}(\Vert u_{0}\Vert_{W^{2,q'}} + \Vert u_{1}\Vert_{q'}), \quad \frac{1}{q} + \frac{1}{q'} = 1, \quad q \in [2,2^{**}],\quad n \geq 1,
\end{equation}
where $u$ is the solution to \eqref{01-p} with $V(x) \equiv 1$ and \eqref{02-p}, and $2^{**} = \infty$ if $n = 1,2,3,4$, and $2^{**} = 2n/(n-4)$ for $n \geq 5$. If, in particular, one takes $q = 2$ in \eqref{p-1}, then one has just boundedness:
\begin{equation}\label{p-2}
\Vert u(t,\cdot)\Vert \leq C(\Vert u_{0}\Vert_{H^{2}} + \Vert u_{1}\Vert)\quad (n \geq 1).
\end{equation}  
So, it seems not trivial to investigate whether $L^{2}$-boundedness holds or not for intermediate potential case of $V(x)$ between $V(x) \equiv 0$ and $V(x) \equiv 1$. In this connection, as for related results about Strichartz estimates of the plate equation, one can cite \cite{CZ} and the references therein.  By a similar multiplier method above one can get the following results. We state them without proof.
\begin{theo}\label{p11} 
Let $n\geq 1$, and $V \in BC({\bf R}^{n})$ satisfies $V(x) > 0$ {\rm (}$\forall x \in {\bf R}^{n}${\rm )}. Suppose that the initial data $[u_{0},u_{1}] \in H^{2}({\bf R}^{n})\times L^{2}({\bf R}^{n})$ further satisfies
\[K_{0,p}^{2} := \int_{{\bf R}^{n}}\frac{\vert u_1(x)\vert^{2}}{V(x)}dx < +\infty.\]
Then the unique solution $u \in C([0,\infty); H^{2}({\bf R}^{n})) \cap C^{1}([0,\infty); L^{2}({\bf R}^{n}))$ to problem \eqref{01-p}-\eqref{02-p} satisfies
\[\Vert u(t,\cdot)\Vert^{2}  \leq C(\Vert u_{0}\Vert^{2} + K_{0,p}^{2}), \quad (t \geq 0)\]
with some generous constant $C > 0$.
\end{theo}

One can also obtain the following result as in Theorem \ref{h12}. The proof is similar, so one states the result without proof. 
\begin{theo}\label{p12} 
Let $\Omega \subset {\bf R}^{n}$ be an open set, and one assumes that $V \in BC({\bf R}^{n})$ satisfies $V(x) > 0$ for $\forall x \in \Omega$, and $V(x) \geq 0$ for $\forall x \in {\bf R}^{n}$. If $u_{1} \in L^{2}({\bf R}^{n})$ additionally satisfies ${\rm supp}u_{1} \subset \Omega$ and
\[J_{0,p}^{2} := \int_{\Omega}\frac{\vert u_{1}(x)\vert^{2}}{V(x)}dx < +\infty,\]
then the unique solution $u \in C([0,\infty); H^{2}({\bf R}^{n})) \cap C^{1}([0,\infty); L^{2}({\bf R}^{n}))$ to problem \eqref{01-p}-\eqref{02-p} satisfies
\[\Vert u(t,\cdot)\Vert^{2}  \leq C(\Vert u_{0}\Vert^{2} + J_{0,p}^{2}), \quad (t \geq 0)\]
with some generous constant $C > 0$.
\end{theo}
\begin{rem}\,{\rm If the potential $V(x)$ is effective only a little, then the growth property occured in \eqref{04-p} can be erased in some cases. Furthermore, the results of Theorems \ref{p11} and \ref{p12} can be viewed as a generalization of \cite[(2.2) with $q=2$ of Theorem 2.1]{L} to the general potential case.} 
\end{rem}
\begin{rem}{\rm By the similar multiplier method one can get the decay estimates of the total energy and $L^{2}$-norm of the solution itself to the damped plate equation:
\[u_{tt} + \Delta^{2}u + V(x)u + a(x)u_{t} = 0, \quad t > 0,\quad x \in {\bf R}^{n},\]
however, this will be left to the readers' exercise.}
\end{rem}

{\bf Notation.}\, We denote the $L^{p}$-norm of $u \in L^{p}({\bf R}^{n})$ by $\Vert u\Vert_{p}$, and in particular, we use $\Vert u\Vert := \Vert u\Vert_{2}$. One sets $B_{R} := \{x \in {\bf R}^{n}\,:\,\vert x\vert < R\}$, and $(f,g) := \displaystyle{\int_{{\bf R}^{n}}f(x)g(x)dx}$ denotes the usual $L^{2}$-inner product of $f,g \in L^{2}({\bf R}^{n})$. We denote $f \in BC({\bf R}^{n})$ $\Leftrightarrow$ $f(x)$ is bounded and continuous in ${\bf R}^{n}$.\\ 

The rest of this paper is organized into three sections. Section 2 is dedicated to prove Theorems \ref{t01} and \ref{t02}, and we prove Theorems \ref{t11} and \ref{t12}. Section 4 is devotd to the proof of Theorems \ref{h11} and \ref{h12}.

\section{Proof of Theorems \ref{t01} and \ref{t02}}

The proofs of a series of Theorems are done by the modified method of \cite{IM}. The method itself of \cite{IM} is a modified version of celebrated paper by \cite{mora}.\\

Let us first prove Theorem \ref{t01}.\\ 

\underline{{\it Proof of Theorem 1.1.}}\, It suffices to assume $[u_{0},u_{1}] \in C_{0}^{\infty}({\bf R}^{n})\times C_{0}^{\infty}({\bf R}^{n})$ to prove results. Then, the corresponding solution becomes sufficiently smooth and vanishs outside of $B_{L+t}$ with some $L > 0$ (cf. \cite{N-2}).

Indeed, as in \cite{IM} set
\[v(t,x) := \int_{0}^{t}u(s,x)ds,\] 
where $u(t,x)$ is the solution to problem (1.1)-(1.2). Then, the function $v(t,x)$ satisfies
\begin{equation}\label{07}
v_{tt}(t,x) - \Delta v(t,x) + V(x)v(t,x) = u_{1}(x), \quad t > 0,\quad x \in {\bf R}^{n},
\end{equation}
\begin{equation}\label{08}
v(0,x)= 0,\hspace{0,50cm}v_{t}(0,x) = u_{0}(x)\quad x \in {\bf R}^{n}.
\end{equation} 
Multiplying the both sides of \eqref{07} by $v_{t}$ and integrating it over $[0,t]\times {\bf R}^{n}$ because of \eqref{08} one has arrived at the idntity such that
\begin{equation}\label{09}
\frac{1}{2}\Vert v_{t}(t,\cdot)\Vert^{2} + \frac{1}{2}\Vert\nabla v(t,\cdot)\Vert^{2} + \frac{1}{2}\int_{{\bf R}^{n}}V(x)\vert v(t,x)\vert^{2}dx = \frac{1}{2}\Vert u_{0}\Vert^{2} + (u_{1}, v(t,\cdot)), \quad t > 0.
\end{equation}
Now, let us estimate the final term of the right hand side of \eqref{09} to absorb it into the left hand side. This part is a crux of our argument. At that case one never relies on the Hardy-type inequality as is usulally used in \cite{IM} in order to include the low dimensional case $n = 1,2$. By the Schwarz inequality one has
\[\vert (u_{1}, v(t,\cdot))\vert \leq \int_{{\bf R}^{n}}\vert u_{1}(x)\vert\vert v(t,x)\vert dx = \int_{{\bf R}^{n}}\frac{\vert u_{1}(x)\vert}{\sqrt{V(x)}}\left(\sqrt{V(x)}\vert\vert v(t,x)\vert\right)dx\]
\[\leq \left(\int_{{\bf R}^{n}}\frac{\vert u_{1}(x)\vert^{2}}{V(x)}dx \right)^{1/2}\left(\int_{{\bf R}^{n}}V(x)\vert v(t,x)\vert^{2}dx\right)^{1/2}\]
\begin{equation}\label{10}
\leq \int_{{\bf R}^{n}}\frac{\vert u_{1}(x)\vert^{2}}{V(x)}dx + \frac{1}{4}\int_{{\bf R}^{n}}V(x)\vert v(t,x)\vert^{2}dx.
\end{equation}
Thus \eqref{09} and \eqref{10} imply the desired estimae:
\begin{equation}\label{11}
\frac{1}{2}\Vert v_{t}(t,\cdot)\Vert^{2} + \frac{1}{2}\Vert\nabla v(t,\cdot)\Vert^{2} + \frac{1}{4}\int_{{\bf R}^{n}}V(x)\vert v(t,x)\vert^{2}dx \leq \frac{1}{2}\Vert u_{0}\Vert^{2} + \int_{{\bf R}^{n}}\frac{\vert u_{1}(x)\vert^{2}}{V(x)}dx,
\end{equation}
because of $v_{t} = u$. It should be repeated that the agument above holds true for all $n \geq 1$ provided that
\[\int_{{\bf R}^{n}}\frac{\vert u_{1}(x)\vert^{2}}{V(x)}dx < +\infty.\]
The help of the potential $V(x)$ is essential.  
\hfill
$\Box$

Let us first prove Theorem \ref{t02}.\\ 

\underline{{\it Proof of Theorem 1.2.}}\,Except for the argument on the part \eqref{10} its proof is similar, so it suffices to compute the corresponding part for \eqref{10}. Indeed, because of the assumption on the initial velocity $u_{1}$ one has
\[\vert (u_{1}, v(t,\cdot))\vert \leq \int_{\Omega}\frac{\vert u_{1}(x)\vert}{\sqrt{V(x)}}\left(\sqrt{V(x)}\vert\vert v(t,x)\vert\right)dx\]
\[\leq \left(\int_{\Omega}\frac{\vert u_{1}(x)\vert^{2}}{V(x)}dx \right)^{1/2}\left(\int_{\Omega}V(x)\vert v(t,x)\vert^{2}dx\right)^{1/2}\]
\begin{equation}\label{12}
\leq \int_{\Omega}\frac{\vert u_{1}(x)\vert^{2}}{V(x)}dx + \frac{1}{4}\int_{{\bf R}^{n}}V(x)\vert v(t,x)\vert^{2}dx.
\end{equation}
\eqref{12} shows the desired estimate.

\hfill
$\Box$


\section{Damped waves}

In this section, we shall prove Theorems \ref{t11} and \ref{t12} by relying on the powerful method coming from \cite{IM}.\\

In this case it is sufficient to suppose $[u_{0},u_{1}] \in C_{0}^{\infty}({\bf R}^{n})\times C_{0}^{\infty}({\bf R}^{n})$. Then, the corresponding solution becomes sufficiently smooth and vanishs outside of $B_{L+t}$ with some $L > 0$.

To begin with, note that the following energy identity holds:
\begin{equation}\label{00-w}
E(t) + \int_{0}^{t}\Vert u_{t}(s,\cdot)\Vert^{2}ds = E(0).
\end{equation}

Multiplying both sides of \eqref{01-w} by $u$ and integrating it on $[0,t]\times {\bf R}^{n}$ one has the identity:
\[\int_{0}^{t}\left(\Vert\nabla u(s,\cdot)\Vert^{2} + \Vert \sqrt{V(\cdot)}u(s,\cdot)\Vert^{2} \right)ds + \frac{1}{2}\Vert u(t,\cdot)\Vert^{2}\]
\begin{equation}\label{03-w}
= \int_{0}^{t}\Vert u_{t}(s,\cdot)\Vert^{2}ds + (u_{1},u_{0}) - (u_{t}(t,\cdot),u(t,\cdot)) + \frac{1}{2}\Vert u_{0}\Vert^{2}\quad (t \geq 0).
\end{equation}
Since 
$$- (u_{t}(t,\cdot),u(t,\cdot)) \leq \Vert u_{t}(t,\cdot)\Vert^{2} + \frac{1}{4}\Vert u(t,\cdot)\Vert^{2}$$ 
from \eqref{03-w} one has
\[\int_{0}^{t}\left(\Vert\nabla u(s,\cdot)\Vert^{2} + \Vert \sqrt{V(\cdot)}u(s,\cdot)\Vert^{2} \right)ds + \frac{1}{4}\Vert u(t,\cdot)\Vert^{2}\]
\begin{equation}\label{04-w}
\leq \int_{0}^{t}\Vert u_{t}(s,\cdot)\Vert^{2}ds + (u_{1},u_{0}) + \frac{1}{2}\Vert u_{0}\Vert^{2} + \Vert u_{t}(t,\cdot)\Vert^{2} \quad (t \geq 0).
\end{equation}
Thus, from \eqref{00-w} and \eqref{04-w} it follows that
\[\int_{0}^{t}\left(\Vert\nabla u(s,\cdot)\Vert^{2} + \Vert \sqrt{V(\cdot)}u(s,\cdot)\Vert^{2} \right)ds + \frac{1}{4}\Vert u(t,\cdot)\Vert^{2} \leq 3E(0) + \frac{1}{2}\Vert u_{0}\Vert^{2} + (u_{1},u_{0}) =: I_{0}^{2},\]
\[\int_{0}^{t}\Vert u_{t}(s,\cdot)\Vert^{2}ds \leq E(0),\]
so that one has arrived at the following result.
\begin{lem}\label{05-w}\,Let $n \geq 1$. Then, it holds that
\[\Vert u(t,\cdot)\Vert \leq 2I_{0}, \quad \int_{0}^{t}E(s)ds \leq I_{0}^{2} + E(0)\quad (t \geq 0),\]
where
\[I_{0} := \sqrt{3E(0) + \frac{1}{2}\Vert u_{0}\Vert^{2} + (u_{1},u_{0})}. \]
\end{lem} 

Next, multiplying both sides of \eqref{01-w} by $(1+t)u_{t}$, and integrating it over $[0,t]\times {\bf R}^{n}$ one has the identity:
\begin{equation}\label{06-w}
(1+t)E(t) + \int_{0}^{t}(1+s)\Vert u_{t}(s,\cdot)\Vert^{2}ds = E(0) + \int_{0}^{t}E(s)ds.
\end{equation}
Furthermore, multiplying both sides of \eqref{01-w} by $(1+t)u$ and integrating it on $[0,t]\times {\bf R}^{n}$ one has the identity:
\[\int_{0}^{t}(1+s)\left(\Vert\nabla u(s,\cdot)\Vert^{2} + \Vert \sqrt{V(\cdot)}u(s,\cdot)\Vert^{2} \right)ds + \frac{(1+t)}{2}\Vert u(t,\cdot)\Vert^{2}\]
\begin{equation}\label{07-w}
= \int_{0}^{t}(1+s)\Vert u_{t}(s,\cdot)\Vert^{2}ds + \frac{1}{2}\Vert u(t,\cdot)\Vert^{2} + (u_{1},u_{0}) - (1+t)(u_{t}(t,\cdot),u(t,\cdot)) + \frac{1}{2}\int_{0}^{t}\Vert u(s,\cdot)\Vert^{2}ds\quad (t \geq 0).
\end{equation}
While, as in the former part one has
$$-(1+t) (u_{t}(t,\cdot),u(t,\cdot)) \leq (1+t)\Vert u_{t}(t,\cdot)\Vert^{2} + \frac{(1+t)}{4}\Vert u(t,\cdot)\Vert^{2}.$$
Therefore, together with \eqref{06-w} and \eqref{07-w} one can get
\[\int_{0}^{t}(1+s)\left(\Vert\nabla u(s,\cdot)\Vert^{2} + \Vert \sqrt{V(\cdot)}u(s,\cdot)\Vert^{2} \right)ds + \frac{(1+t)}{4}\Vert u(t,\cdot)\Vert^{2}\]
\[\leq \int_{0}^{t}(1+s)\Vert u_{t}(s,\cdot)\Vert^{2}ds + \frac{1}{2}\Vert u(t,\cdot)\Vert^{2} + (u_{1},u_{0}) + (1+t)\Vert u_{t}(t,\cdot)\Vert^{2} + \int_{0}^{t}\Vert u(s,\cdot)\Vert^{2}ds\]
\[\leq (E(0) + \int_{0}^{t}E(s)ds) + \frac{1}{2}\Vert u(t,\cdot)\Vert^{2} + (u_{1},u_{0}) + 2(E(0) + \int_{0}^{t}E(s)ds) + \int_{0}^{t}\Vert u(s,\cdot)\Vert^{2}ds\]
\begin{equation}\label{08-w}
\leq 3E(0) + 3\int_{0}^{t}E(s)ds + (u_{1},u_{0}) + \frac{1}{2}\Vert u(t,\cdot)\Vert^{2} + \int_{0}^{t}\Vert u(s,\cdot)\Vert^{2}ds \quad (t \geq 0).
\end{equation}
\eqref{06-w} and \eqref{08-w} implies the following lemma.
\begin{lem}\label{09-w}\,Let $n \geq 1$. Then, it holds that
\[\int_{0}^{t}(1+s)\left(\Vert\nabla u(s,\cdot)\Vert^{2} + \Vert \sqrt{V(\cdot)}u(s,\cdot)\Vert^{2} \right)ds + \frac{(1+t)}{4}\Vert u(t,\cdot)\Vert^{2}\]
\[\leq 3E(0) + 3\int_{0}^{t}E(s)ds + (u_{1},u_{0}) + \frac{1}{2}\Vert u(t,\cdot)\Vert^{2} + \int_{0}^{t}\Vert u(s,\cdot)\Vert^{2}ds \quad (t \geq 0),\]
and
\[\int_{0}^{t}(1+s)\Vert u_{t}(s,\cdot)\Vert^{2}ds \leq E(0) + \int_{0}^{t}E(s)ds \quad (t \geq 0).\]
\end{lem} 
In order to get the fast decay estimate one has to note the following inequality:
\[\frac{d}{dt}\left((1+t)^{2}E(t) \right) \leq 2(1+t)E(t),\]
so that from Lemmas \ref{05-w} and \ref{09-w} one has
\[(1+t)^{2}E(t) \leq 5E(0) + 4\int_{0}^{t}E(s)ds + (u_{1},u_{0}) + \frac{1}{2}\Vert u(t,\cdot)\Vert^{2} + \int_{0}^{t}\Vert u(s,\cdot)\Vert^{2}ds,\]
so the lemma below holds.
\begin{lem}\label{10-w}\,Let $n \geq 1$. Then, it holds that
\[(1+t)^{2}E(t) \leq 9E(0) + 6I_{0}^{2}  + (u_{1},u_{0}) + \int_{0}^{t}\Vert u(s,\cdot)\Vert^{2}ds \quad (t \geq 0),\]
\[\frac{(1+t)}{4}\Vert u(t,\cdot)\Vert^{2} \leq 6E(0) + 5I_{0}^{2}  + (u_{1},u_{0}) + \int_{0}^{t}\Vert u(s,\cdot)\Vert^{2}ds \quad (t \geq 0).\]
\end{lem} 
From the observation on Lemma \ref{10-w} one can notice that to get the fast energy and $L^{2}$-decay it suffices to obtain the a priori bound for the quantity 
\[\int_{0}^{t}\Vert u(s,\cdot)\Vert^{2}ds.\]
This is a crux of this paper. For this end, the existence of potentoal function $V(x)$ plays an essential role, in particular, on the treatment for the low dimensional case $n = 1,2$ because we never relies on the Hardy type inequalty as is previously discussed in \cite{IM}. Additionally, the assumption given to the potential is much generous than known results as in \cite{N}. 
\begin{lem}\label{11-w}\,Let $n \geq 1$. Then, it holds that
\[\int_{0}^{t}\Vert u(s,\cdot)\Vert^{2}ds \leq  C\left(\Vert u_{0}\Vert^{2} + \int_{{\bf R}^{n}}\frac{\vert u_{1}(x)+u_{0}(x)\vert^{2}}{V(x)}dx\right)\quad (t \geq 0),\]
where $C > 0$ is a generous constant.
\end{lem} 

\underline{{\it Proof of Lemma \ref{11-w}.}}\, The proof is almost similar to the non-damped case developed in Theorem 1.1. In fact, set 
\[v(t,x) := \int_{0}^{t}u(s,x)ds,\] 
where $u(t,x)$ is the solution to problem \eqref{01-w}-\eqref{02-w}. Then, the function $v(t,x)$ satisfies
\begin{equation}\label{12-w}
v_{tt}(t,x) - \Delta v(t,x) + V(x)v(t,x) + v_{t}(t,x) = u_{1}(x) + u_{0}(x), \quad t > 0,\quad x \in {\bf R}^{n},
\end{equation}
\begin{equation}\label{13-w}
v(0,x)= 0,\hspace{0,50cm}v_{t}(0,x) = u_{0}(x)\quad x \in {\bf R}^{n}.
\end{equation} 
Multiplying the both sides of \eqref{12-w} by $v_{t}$ and integrating it over $[0,t]\times {\bf R}^{n}$ because of \eqref{13-w} one has arrived at the identity such that

\[\frac{1}{2}\Vert v_{t}(t,\cdot)\Vert^{2} + \frac{1}{2}\Vert\nabla v(t,\cdot)\Vert^{2} + \frac{1}{2}\int_{{\bf R}^{n}}V(x)\vert v(t,x)\vert^{2}dx + \int_{0}^{t}\Vert v_{s}(s,\cdot)\Vert^{2}ds\]
\begin{equation}\label{14-w}
= \frac{1}{2}\Vert u_{0}\Vert^{2} + (u_{1} + u_{0}, v(t,\cdot)), \quad t \geq 0.
\end{equation}
Now, let us estimate the final term of the right hand side of \eqref{14-w} to absorb it into the left hand side. By the Schwarz inequality one has
\[\vert (u_{1}+u_{0}, v(t,\cdot))\vert \leq \int_{{\bf R}^{n}}\vert u_{1}(x)+u_{0}(x)\vert\vert v(t,x)\vert dx = \int_{{\bf R}^{n}}\frac{\vert u_{1}(x)+u_{0}(x)\vert}{\sqrt{V(x)}}\left(\sqrt{V(x)}\vert\vert v(t,x)\vert\right)dx\]
\[\leq \left(\int_{{\bf R}^{n}}\frac{\vert u_{1}(x)+u_{0}(x)\vert^{2}}{V(x)}dx \right)^{1/2}\left(\int_{{\bf R}^{n}}V(x)\vert v(t,x)\vert^{2}dx\right)^{1/2}\]
\begin{equation}\label{15-w}
\leq \int_{{\bf R}^{n}}\frac{\vert u_{1}(x)+u_{0}(x)\vert^{2}}{V(x)}dx + \frac{1}{4}\int_{{\bf R}^{n}}V(x)\vert v(t,x)\vert^{2}dx.
\end{equation}
Thus \eqref{14-w} and \eqref{15-w} imply the desired estimae:

\[\frac{1}{2}\Vert v_{t}(t,\cdot)\Vert^{2} + \frac{1}{2}\Vert\nabla v(t,\cdot)\Vert^{2} + \frac{1}{4}\int_{{\bf R}^{n}}V(x)\vert v(t,x)\vert^{2}dx + \int_{0}^{t}\Vert v_{s}(s,\cdot)\Vert^{2}ds\]
\begin{equation}\label{16-w}
\leq \frac{1}{2}\Vert u_{0}\Vert^{2} + \int_{{\bf R}^{n}}\frac{\vert u_{1}(x)+u_{0}(x)\vert^{2}}{V(x)}dx,
\end{equation}
because of $v_{t} = u$. It should be repeated that the agument above holds true for all $n \geq 1$ provided that
\[\int_{{\bf R}^{n}}\frac{\vert u_{1}(x)+u_{0}(x)\vert^{2}}{V(x)}dx < +\infty.\]
The help of the potential $V(x)$ is essential.  
\hfill
$\Box$\\

Proof of Theorem \ref{t11} is a direct consequence of Lemmas \ref{10-w} and \ref{11-w}. Furthermore, one can prove similarly Theorem \ref{t12} by modifying the proof of Theorem \ref{t11} above basing on the proof of Theorem \ref{02}.


\section{Heat equations}

In this section, we shall prove Theorems \ref{h11} and \ref{h12} similarly to Section 3.\\
For the moment it is sufficient to suppose $u_{0} \in C_{0}^{\infty}({\bf R}^{n})$. Then the corresponding solution $u(t,x)$ is smoother to guarantee the integration by parts below.\\

First, multiplying both sides of \eqref{01-h} by $u$, and integrating over ${\bf R}^{n}$ one has
\[\frac{d}{dt}\Vert u(t,\cdot)\Vert^{2} + \Vert\nabla u(t,\cdot)\Vert^{2} + \int_{{\bf R}^{n}}V(x)\vert u(t,x)\vert^{2}dx = 0,\quad t > 0,\]
so that
\[
\frac{d}{dt}\Vert u(t,\cdot)\Vert^{2} \leq 0,\quad t > 0.
\]
Thus, it follows that
\[\frac{d}{dt}\left((1+t)\Vert u(t,\cdot)\Vert^{2}\right) \leq \Vert u(t,\cdot)\Vert^{2}.\] 
By integrating it over $[0,t]$ one has 
\begin{equation}\label{h-10}
(1+t)\Vert u(t,\cdot)\Vert^{2} \leq \Vert u_{0}\Vert^{2} + \int_{0}^{t}\Vert u(s,\cdot)\Vert^{2}ds, \quad t > 0.
\end{equation}
While, as in \cite{IM} one sets again
\[v(t,x) := \int_{0}^{t}u(s,x)ds.\]
Then, the function $v(t,x)$ satisfies
\begin{equation}\label{h-11}
v_{t}- \Delta v + V(x)v = u_{0}, \quad t > 0,\quad x \in {\bf R}^{n},
\end{equation}
\begin{equation}\label{h-12}
v(0,x)=0,\quad x \in {\bf R}^{n}.
\end{equation} 
Multiplying both sides of \eqref{h-11} by $v_{t}$ and integrating it over $[0, t]\times {\bf R}^{n}$ because of \eqref{h-12} one has
\[\int_{0}^{t}\Vert v_{s}(s,\cdot)\Vert^{2}ds + \frac{1}{2}\Vert \nabla v(t,\cdot)\Vert^{2} + \frac{1}{2}\Vert \sqrt{V}v(t,\cdot)\Vert^{2}\]
\begin{equation}\label{h-13-0}
= (u_{0}, v(t,\cdot)).
\end{equation}
Bcause of $v_{s} = u$, \eqref{h-13-0} implies
\[\int_{0}^{t}\Vert u(s,\cdot)\Vert^{2}ds + \frac{1}{2}\Vert \sqrt{V}v(t,\cdot)\Vert^{2}\]
\begin{equation}\label{h-13}
\leq (u_{0}, v(t,\cdot)) \quad (t > 0).
\end{equation}
Note that \eqref{h-10} and \eqref{h-13} can be true in the case of $u_{0} \in L^{2}({\bf R}^{n})$ and the corresponding solution $u$ (and $v$) by density argument (cf. \cite[Theorem VII.7]{B}).\\
Now we apply the argument as in the proof of Section 3. Then, one can get
\begin{equation}\label{h-14}
(u_{0}, v(t,\cdot)) \leq \int_{{\bf R}^{n}}\frac{\vert u_{0}(x)\vert^{2}}{V(x)}dx + \frac{1}{4}\int_{{\bf R}^{n}}V(x)\vert v(t,x)\vert^{2}dx.
\end{equation} 
\eqref{h-13} and \eqref{h-14} imply
\[\int_{0}^{t}\Vert u(s,\cdot)\Vert^{2}ds + \frac{1}{4}\Vert \sqrt{V}v(t,\cdot)\Vert^{2} \leq \int_{{\bf R}^{n}}\frac{\vert u_{0}(x)\vert^{2}}{V(x)}dx,\quad t > 0.\]
Thus, one has arrived at the important estimate
\begin{equation}\label{h-15}
\int_{0}^{t}\Vert u(s,\cdot)\Vert^{2}ds \leq \int_{{\bf R}^{n}}\frac{\vert u_{0}(x)\vert^{2}}{V(x)}dx < +\infty, \quad t > 0.
\end{equation} 
The desired statement of Theorem \ref{h11} can be derived by \eqref{h-10} and \eqref{h-15}. \\

The proof of Theorem \ref{h12} can be given similarly by using the argument developed in the proof of Theorem \ref{t02}.
\hfill
$\Box$
\par
\vspace{0.5cm}
\noindent{\em Acknowledgement.}
\smallskip
The work of the author was supported in part by Grant-in-Aid for Scientific Research (C) 20K03682  of JSPS. 


\end{document}